\documentclass{article}

\title{\LARGE \textbf{Long Cycles in 1-tough Graphs}}
\author{Zh.G. Nikoghosyan\footnote{G.G. Nicoghossian (up to 1997)}}

\begin{document}

\maketitle

\begin{abstract}
In 1952, Dirac proved that every 2-connected graph with minimum degree $\delta$ either is hamiltonian or contains a cycle of length at least $2\delta$. In 1986, Bauer and Schmeichel enlarged the bound $2\delta$ to $2\delta+2$ under additional 1-tough condition - an alternative and more natural necessary condition for a graph to be hamiltonian.  In fact, the bound $2\delta+2$ is sharp for a graph on $n$ vertices when $n\equiv 1(mod\ 3)$. In this paper we present the final version of this result which is sharp for each $n$: every 1-tough graph either is hamiltonian or contains a cycle of length at least $2\delta+2$ when $n\equiv 1(mod\ 3)$, at least $2\delta+3$ when $n\equiv 2(mod\ 3)$ or $n\equiv 1(mod\ 4)$, and at least $2\delta+4$ otherwise.  \\


\end{abstract}

\section{Introduction}

Only finite undirected graphs without loops or multiple edges are considered. A good reference for any undefined terms is \cite{[2]}. 

The earliest two theoretical results in hamiltonian graph theory were developed in 1952 due to Dirac \cite{[4]} in forms of a sufficient condition for a graph to be hamiltonian and a lower bound for the circumference $c$ (the length of a longest cycle of a graph), respectively, based on order $n$ and minimum degree $\delta$.  \\

\noindent\textbf{Theorem A \cite{[4]}.} Every graph with $\delta\ge n/2$ is hamiltonian. \\

\noindent\textbf{Theorem B \cite{[4]}.} Every 2-connected graph either is hamiltonian or $c\ge 2\delta$. \\

In 1973, Chv\'{a}tal  \cite{[3]} introduced the concept of toughness. Since then a lot of research has been done towards finding the exact analogs of classical hamiltonian results under additional 1-tough condition - an alternative and more natural  necessary condition for a graph to be hamiltonian. 

In 1978, Jung \cite{[5]} established the analog of Theorem A for 1-tough graphs.\\

\noindent\textbf{Theorem C \cite{[5]}.} Every 1-tough graph on $n\ge11$ vertices with $\delta\ge(n-4)/2$ is hamiltonian.\\

In 1986, Bauer and Schmeichel \cite{[1]} proved that the bound $2\delta$ in Theorem B can be enlarged to $2\delta+2$ for 1-tough graphs.\\

\noindent\textbf{Theorem D \cite{[1]}.} Every 1-tough graph either is hamiltonian or $c\ge 2\delta+2$.   \\

In fact, Theorem D is sharp when $n\equiv 1(mod\ 3)$. In this paper we present the final version of Theorem D which is sharp for each $n$.\\

\noindent\textbf{Theorem 1}. Every 1-tough graph either is hamiltonian or
$$
c\ge\left\{ 
\begin{array}{lll}
2\delta+2 & \mbox{when} &  n\equiv 1(mod\ 3), \\ 2\delta+3 & \mbox{when} & %
n\equiv 2(mod\ 3) \ \ \mbox{or} \ \ n\equiv 1(mod\ 4), \\ 2\delta+4 & \mbox{ } & \mbox{otherwise}. 
\end{array}
\right. 
$$
Theorem 1 is sharp for each $n$. To see this, let $H_1,H_2,...,H_h$ be disjoint complete graphs with distinct vertices $x_i,y_i\in V(H_i)$ $(i=1,2,...,h)$. Form a new graph $H(t_1,t_2,...,t_h)$ by identifying the vertices $x_1,x_2,...,x_h$ and adding all possible edges between $y_1,y_2,...,y_h$, where $t_i=|V(H_i)|$ $(i=1,2,...,h)$. The graph $H(\delta-1,\delta-1,\delta-1)$ shows that the bound $2\delta+2$ in Theorem 1 cannot be replaced by $2\delta+3$ when $n\equiv 1(mod\ 3)$. Next, the graphs $H(\delta,\delta-1,\delta-1)$ and $H(\delta-1,\delta-1,\delta-1,\delta-1)$ show that the bound $2\delta+3$ cannot be replaced by $2\delta+4$ when $n\equiv 2(mod\ 3)$ or $n\equiv 1(mod\ 4)$. Finally, the graph $H(\delta,\delta,\delta-1)$ shows that the bound $2\delta+4$ cannot be replaced by $2\delta+5$.

\section{Notations and preliminaries}

The set of vertices of a graph $G$ is denoted by $V(G)$ and the set of edges by $E(G)$. For $S$ a subset of $V(G)$, we denote by $G\backslash S$ the maximum subgraph of $G$ with vertex set $V(G)\backslash S$. We write $\langle S\rangle$ for the subgraph of $G$ induced by $S$. For a subgraph $H$ of $G$ we use $G\backslash H$ short for $G\backslash V(H)$. The neighborhood and the degree of a vertex $x\in V(G)$ will be denoted by $N(x)$ and $d(x)$, respectively.  Furthermore, for a subgraph $H$ of $G$ and $x\in V(G)$, we define $N_H(x)=N(x)\cap V(H)$ and  $d_H(x)=|N_H(x)|$. Let $s(G)$ denote the number of components of a graph $G$. A graph $G$ is 1-tough if $|S|\ge s(G\backslash S)$ for every subset $S$ of the vertex set $V(G)$ with $s(G\backslash S)>1$.  A graph $G$ on $n$ vertices is hamiltonian if $G$ contains a Hamilton cycle, i.e. a cycle of length $n$.

Paths and cycles in a graph $G$ are considered as subgraphs of $G$. If $Q$ is a path or a cycle, then the length of $Q$, denoted by $|Q|$, is $|E(Q)|$. We write $Q$ with a given orientation by $\overrightarrow{Q}$. For $x,y\in V(Q)$, we denote by $x\overrightarrow{Q}y$ the subpath of $Q$ in the chosen direction from $x$ to $y$. For $x\in V(C)$, we denote the $h$-th successor and the $h$-th predecessor of $x$ on $\overrightarrow{C}$ by $x^{+h}$ and $x^{-h}$, respectively. We abbreviate $x^{+1}$ and $x^{-1}$ by $x^+$ and $x^-$, respectively.  For each $X\subset V(C)$, we define $X^{+}=\{x^{+}|x\in X\}$ and $X^{-}=\{x^{-}|x\in X\}$. \\

\noindent\textbf{Special definitions}. Let $G$ be a graph, $C$ a longest cycle in $G$ and $P=x\overrightarrow{P}y$ a longest path in $G\backslash C$ of length $\overline{p}\ge0$. Let $\xi_1,\xi_2,...,\xi_s$ be the elements of $N_C(x)\cup N_C(y)$ occuring on $\overrightarrow{C}$ in a consecutive order. Set 
$$
I_i=\xi_i\overrightarrow{C}\xi_{i+1},       \     I_i^\ast=\xi_i^+\overrightarrow{C}\xi_{i+1}^-   \   \  (i=1,2,...,s),
$$
where $\xi_{s+1}=\xi_1$. 

$(1)$ The segments  $I_1,I_2,...,I_s$ are called elementary segments on $C$ induced by $N_C(x)\cup N_C(y)$.

$(2)$ We call a path $L=z\overrightarrow{L}w$ an intermediate path between two distinct elementary segments $I_a$ and $I_b$, if
$$
z\in V(I_a^\ast),   \    w\in V(I_b^\ast),    \    V(L)\cap V(C\cup P)=\{z,w\}.
$$

$(3)$ Define $\Upsilon(I_{i_1},I_{i_2},...,I_{i_t})$ to be the set of all intermediate paths between elementary segments  $I_{i_1},I_{i_2},...,I_{i_t}$.\\

$(4)$ If $\Upsilon(I_1,...,I_s)\subseteq E$ then the maximum number of intermediate independent edges (having no a common vertex) in $\Upsilon(I_1,...,I_s)$ will be denoted by $\mu(\Upsilon)$.\\

$(5)$ We say that two intermediate independent edges $w_1w_2, w_3w_4$ have a  crossing, if either $w_1,w_3,w_2,w_4$ or $w_1,w_4,w_2,w_3$ occur on $\overrightarrow{C}$ in a consecutive order.\\

\noindent\textbf{Lemma 1.} Let $G$ be a graph, $C$ a longest cycle in $G$ and $P=x\overrightarrow{P}y$ a longest path in $G\backslash C$ of length $\overline{p}\ge1$. If  $|N_C(x)|\ge2$, $|N_C(y)|\ge2$ and  $N_C(x)\not=N_C(y)$ then
$$
c\ge\left\{ 
\begin{array}{lll}
3\delta+\max\{\sigma_1, \sigma_2\}-1\ge3\delta & \mbox{if} & \mbox{ }\overline{p}=1, \\ 4\delta-2\overline{p} & \mbox{if} & \mbox{ }%
\overline{p}\ge2, 
\end{array}
\right. 
$$
where $\sigma_1=|N_C(x)\backslash N_C(y)|$ and $\sigma_2=|N_C(y)\backslash N_C(x)|$.\\

\noindent\textbf{Lemma 2.} Let $G$ be a graph, $C$ a longest cycle in $G$ and $P=x\overrightarrow{P}y$ a longest path in $G\backslash C$ of length $\overline{p}\ge0$.  Let $N_C(x)=N_C(y)$, $|N_C(x)|\ge2$ and $f,g\in\{1,...,s\}$.\\

(a1) If $L\in\Upsilon(I_f,I_g)$ then  
$$
|I_f|+|I_g|\ge2\overline{p}+2|L|+4.
$$

(a2) If $\Upsilon(I_f,I_g)\subseteq E(G)$ and $|\Upsilon(I_f,I_g)|=\varepsilon$\ \  for some $\varepsilon\in\{1,2,3\}$ then 
$$
|I_f|+|I_g|\ge2\overline{p}+\varepsilon+5,
$$

(a3) If $\Upsilon(I_f,I_g)\subseteq E(G)$ and $\Upsilon(I_f,I_g)$ contains two independent intermediate edges then 
$$
|I_f|+|I_g|\ge2\overline{p}+8.
$$
\\
The following result is due to Voss \cite{[6]}.  \\

\noindent\textbf{Lemma 3 \cite{[6]}}. Let $G$ be a hamiltonian graph, $\{v_1,v_2,...,v_t\}\subseteq V(G)$ and $d(v_i)\ge t$ $(i=1,2,...,t)$. Then each pair $x,y$ of vertices of $G$ is connected in $G$ by a path of length at least $t$.\\

\section{Proofs}

\noindent\textbf{Proof of Lemma 1}. Put
$$
A_1=N_C(x)\backslash N_C(y),    \    A_2=N_C(y)\backslash N_C(x), \   M=N_C(x)\cap N_C(y).
$$
By the hypothesis,  $N_C(x)\not=N_C(y)$, implying that 
$$
\max \{|A_1|,|A_2|\}\ge1.
$$ 
Let $\xi_1,\xi_2,...,\xi_s$ be the elements of $N_C(x)\cup N_C(y)$ occuring on \overrightarrow{C} in a consecutive order. Put $I_i=\xi_i\overrightarrow{C}\xi_{i+1}$ $(i=1,2,...,s)$, where $\xi_{s+1}=\xi_1$. Clearly, $s=|A_1|+|A_2|+|M|$. Since $C$ is extreme, we have $|I_i|\ge2$ $(i=1,2,...,s)$. Next, if $\{\xi_i,\xi_{i+1}\}\cap M\not=\emptyset$ for some $i\in\{1,2,...,s\}$ then $|I_i|\ge\overline{p}+2$. Further, if either $\xi_i\in A_1$, $\xi_{i+1}\in A_2$ or $\xi_i\in A_2$, $\xi_{i+1}\in A_1$ then again $|I_i|\ge\overline{p}+2$. \\

\textbf{Case 1}. $\overline{p}=1$.

\textbf{Case 1.1}. $|A_i|\ge1$ $(i=1,2)$.

It follows that among  $I_1,I_2,...,I_s$ there are  $|M|+2$ segments of length at least $\overline{p}+2$. Observing also that each of the remaining $s-(|M|+2)$ segments has a length at least 2, we have 
$$
c\ge(\overline{p}+2)(|M|+2)+2(s-|M|-2)
$$
$$
=3(|M|+2)+2(|A_1|+|A_2|-2)=2|A_1|+2|A_2|+3|M|+2.
$$

Since $|A_1|=d(x)-|M|-1$ and  $|A_2|=d(y)-|M|-1$, we have
$$
c\ge2d(x)+2d(y)-|M|-2\ge3\delta+d(x)-|M|-2.
$$
Recalling that $d(x)=|M|+|A_1|+1$, we get 
$$
c\ge3\delta+|A_1|-1=3\delta+\sigma_1-1.
$$
Analogously, $c\ge3\delta+\sigma_2-1$. So, 
$$
c\ge3\delta+\max \{\sigma_1,\sigma_2\}-1\ge3\delta.
$$

\textbf{Case 1.2}. Either $|A_1|\ge1, |A_2|=0$ or $|A_1|=0, |A_2|\ge1$.

Assume w.l.o.g. that $|A_1|\ge1$ and $|A_2|=0$, i.e. $|N_C(y)|=|M|\ge2$ and $s=|A_1|+|M|$ . Hence, among $I_1,I_2,...,I_s$ there are $|M|+1$ segments of length at least $\overline{p}+2=3$. Taking into account that $|M|+1=d(y)$ and each of the remaining $s-(|M|+1)$ segments has a length at least 2, we get
$$
c\ge 3(|M|+1)+2(s-|M|-1)=3d(y)+2(|A_1|-1)
$$
$$
\ge3\delta+|A_1|-1=3\delta+\max\{\sigma_1,\sigma_2\}-1\ge3\delta.
$$

\textbf{Case 2}. $\overline{p}\ge2$.

\textbf{Case 2.1}. $|A_i|\ge1$ $(i=1,2)$.

It follows that among  $I_1,I_2,...,I_s$ there are $|M|+2$ segments of length at least $\overline{p}+2$. Further, since each of the remaining $s-(|M|+2)$ segments has a length at least 2, we get 
$$
c\ge (\overline{p}+2)(|M|+2)+2(s-|M|-2)
$$
$$
=(\overline{p}-2)|M|+(2\overline{p}+4|M|+4)+2(|A_1|+|A_2|-2)
$$
$$
\ge2|A_1|+2|A_2|+4|M|+2\overline{p}.
$$
Observing also that 
$$
|A_1|+|M|+\overline{p}\ge d(x),  \quad   |A_2|+|M|+\overline{p}\ge d(y),  
$$
we have 

$$
2|A_1|+2|A_2|+4|M|+2\overline{p}
$$
$$
\ge 2d(x)+2d(y)-2\overline{p}\ge4\delta-2\overline{p},
$$
implying that $c\ge4\delta-2\overline{p}$.\\

 \textbf{Case 2.2}. Either $|A_1|\ge1, |A_2|=0$ or $|A_1|=0, |A_2|\ge1$.

Assume w.l.o.g. that $|A_1|\ge1$ and $|A_2|=0$, that is  $|N_C(y)|=|M|\ge2$ and $s=|A_1|+|M|$. It follows that among  $I_1,I_2,...,I_s$ there are  $|M|+1$ segments of length at least $\overline{p}+2$. Observing also that $|M|+\overline{p}\ge d(y)\ge\delta$, i.e. $2\overline{p}+4|M|\ge 4\delta-2\overline{p}$, we get 
$$
c\ge(\overline{p}+2)(|M|+1)\ge(\overline{p}-2)(|M|-1)+2\overline{p}+4|M|
$$
$$
\ge 2\overline{p}+4|M|\ge4\delta-2\overline{p}.       \quad            \quad                    \rule{7pt}{6pt} 
$$

\noindent\textbf{Proof of Lemma 2}. Let $\xi_1,\xi_2,...,\xi_s$ be the elements of $N_C(x)$ occuring on $\overrightarrow{C}$ in a consecutive order. Put $I_i=\xi_i\overrightarrow{C}\xi_{i+1}$ $(i=1,2,...,s)$, where $\xi_{s+1}=\xi_1.$  To prove $(a1)$, let $L\in \Upsilon(I_f,I_g)$. Further, let $L=z\overrightarrow{L}w$ with $z\in V(I_f^\ast)$ and   $w\in V(I_g^\ast)$.
Put
$$
|\xi_f\overrightarrow{C}z|=d_1,  \   |z\overrightarrow{C}\xi_{f+1}|=d_2,      \      |\xi_g\overrightarrow{C}w|=d_3,       \     |w\overrightarrow{C}\xi_{g+1}|=d_4,
$$
$$
C^\prime=\xi_fx\overrightarrow{P}y\xi_g\overleftarrow{C}z\overrightarrow{L}w\overrightarrow{C}\xi_f.
$$
Clearly, 
$$
|C^\prime|=|C|-d_1-d_3+|L|+|P|+2.           
$$
Since $C$ is extreme, we have $|C|\ge|C^\prime|$, implying that $d_1+d_3\ge\overline{p}+|L|+2$.  By a symmetric argument, $d_2+d_4\ge\overline{p}+|L|+2$. Hence 
$$
|I_f|+|I_g|=\sum_{i=1}^4d_i\ge2\overline{p}+2|L|+4.
$$

The proof of $(a1)$ is complete. To prove  $(a2)$ and $(a3)$, let  $\Upsilon(I_f,I_g)\subseteq E(G)$ and $|\Upsilon(I_f,I_g)|=\varepsilon$ for some $\varepsilon\in \{1,2,3\}$.\\

\textbf{Case 1}. $\varepsilon=1$.

Let $L\in\Upsilon(I_f,I_g)$, where $|L|=1$. By (a1), 
$$
|I_f|+|I_g|\ge2\overline{p}+2|L|+4=2\overline{p}+6.
$$

\textbf{Case 2}. $\varepsilon=2$.

It follows that $\Upsilon(I_f,I_g)$ consists of two edges $e_1,e_2$. Put $e_1=z_1w_1$ and $e_2=z_2w_2$, where $\{z_1,z_2\}\subseteq V(I_f^\ast)$ and $\{w_1,w_2\}\subseteq V(I_g^\ast)$.\\

\textbf{Case 2.1}. $z_1\not=z_2$ and $w_1\not=w_2$.

Assume w.l.o.g. that $z_1$ and $z_2$ occur in this order on $I_f$. \\

\textbf{Case 2.1.1}. $w_2$ and  $w_1$ occur in this order on $I_g$.

Put
$$
|\xi_f\overrightarrow{C}z_1|=d_1,  \   |z_1\overrightarrow{C}z_2|=d_2,      \    |z_2\overrightarrow{C}\xi_{f+1}|=d_3, 
$$
$$
|\xi_g\overrightarrow{C}w_2|=d_4,       \       |w_2\overrightarrow{C}w_1|=d_5,         \        |w_1\overrightarrow{C}\xi_{g+1}|=d_6, 
$$
$$
C^{\prime}=\xi_f\overrightarrow{C}z_1w_1\overleftarrow{C}w_2z_2\overrightarrow{C}\xi_g x\overrightarrow{P}y\xi_{g+1}\overrightarrow{C}\xi_f.
$$
Clearly, 
$$
|C^{\prime}|=|C|-d_2-d_4-d_6+|\{e_1\}|+|\{e_2\}|+|P|+2
$$
$$
=|C|-d_2-d_4-d_6+\overline{p}+4.
$$
Since $C$ is extreme, we have $|C|\ge |C^{\prime}|$, implying that $d_2+d_4+d_6\ge \overline{p}+4$. By a symmetric argument, $d_1+d_3+d_5\ge\overline{p}+4$. Hence
$$
|I_f|+|I_g|= \sum_{i=1}^6d_i\ge2\overline{p}+8.
$$

\textbf{Case 2.1.2}. $w_1$ and $w_2$ occur in this order on $I_g$.

Putting
$$
C^{\prime}=\xi_f\overrightarrow{C}z_1w_1\overrightarrow{C}w_2z_2\overrightarrow{C}\xi_g x\overrightarrow{P}y\xi_{g+1}\overrightarrow{C}\xi_f,
$$
we can argue as in Case 2.1.1. \\

\textbf{Case 2.2}. Either $z_1=z_2$, $w_1\not=w_2$ or $z_1\not=z_2$, $w_1=w_2$.

Assume w.l.o.g. that $z_1\not=z_2$, $w_1=w_2$ and $z_1, z_2$ occur  in this order on $I_f$. Put
$$
|\xi_f\overrightarrow{C}z_1|=d_1,  \   |z_1\overrightarrow{C}z_2|=d_2,      \    |z_2\overrightarrow{C}\xi_{f+1}|=d_3, 
$$
$$
|\xi_g\overrightarrow{C}w_1|=d_4,       \           |w_1\overrightarrow{C}\xi_{g+1}|=d_5, 
$$
$$
C^{\prime}=\xi_f x\overrightarrow{P}y\xi_g\overleftarrow{C}z_1w_1\overrightarrow{C}\xi_f,
$$
$$
C^{\prime\prime}=\xi_f\overrightarrow{C}z_2w_1\overleftarrow{C}\xi_{f+1}x\overrightarrow{P}y\xi_{g+1}\overrightarrow{C}\xi_f.
$$

Clearly, 
$$
|C^{\prime}|=|C|-d_1-d_4+|\{e_1\}|+|P|+2=|C|-d_1-d_4+\overline{p}+3,
$$
$$
|C^{\prime\prime}|=|C|-d_3-d_5+|\{e_2\}|+|P|+2=|C|-d_3-d_5+\overline{p}+3.
$$
Since $C$ is extreme, $|C|\ge |C^{\prime}|$ and $|C|\ge |C^{\prime\prime}|$, implying that 
$$
d_1+d_4\ge \overline{p}+3,  \    d_3+d_5\ge \overline{p}+3. 
$$
Hence, 
$$
|I_f|+|I_g|= \sum_{i=1}^5d_i\ge d_1+d_3+d_4+d_5+1\ge2\overline{p}+7.
$$

\textbf{Case 3}. $\varepsilon=3$.

It follows that $\Upsilon(I_f,I_g)$ consists of three edges $e_1,e_2,e_3$. Let $e_i=z_iw_i$ $(i=1,2,3)$, where $\{z_1,z_2,z_3\}\subseteq V(I_f^\ast)$ and $\{w_1,w_2,w_3\}\subseteq V(I_g^\ast)$. If there are two independent edges among $e_1,e_2,e_3$ then we can argue as in Case 2.1. Otherwise, we can assume w.l.o.g. that $w_1=w_2=w_3$ and $z_1,z_2,z_3$ occur in this order on $I_f$. Put
$$
|\xi_f\overrightarrow{C}z_1|=d_1,  \   |z_1\overrightarrow{C}z_2|=d_2,      \    |z_2\overrightarrow{C}z_3|=d_3,      
$$
$$
|z_3\overrightarrow{C}\xi_{f+1}|=d_4,      \    |\xi_g\overrightarrow{C}w_1|=d_5,       \           |w_1\overrightarrow{C}\xi_{g+1}|=d_6, 
$$
$$
C^{\prime}=\xi_f x\overrightarrow{P}y\xi_g\overleftarrow{C}z_1w_1\overrightarrow{C}\xi_f,
$$
$$
C^{\prime\prime}=\xi_f\overrightarrow{C}z_3w_1\overleftarrow{C}\xi_{f+1}x\overrightarrow{P}y\xi_{g+1}\overrightarrow{C}\xi_f.
$$

Clearly, 
$$
|C^{\prime}|=|C|-d_1-d_5+|\{e_1\}|+\overline{p}+2,
$$
$$
|C^{\prime\prime}|=|C|-d_4-d_6+|\{e_3\}|+\overline{p}+2.
$$
Since $C$ is extreme, we have $|C|\ge |C^{\prime}|$ and $|C|\ge |C^{\prime\prime}|$, implying that 
$$
d_1+d_5\ge \overline{p}+3,   \    d_4+d_6\ge \overline{p}+3. 
$$
Hence, 
$$
|I_f|+|I_g|= \sum_{i=1}^6d_i\ge d_1+d_4+d_5+d_6+2\ge2\overline{p}+8.     \quad      \quad    \rule{7pt}{6pt} 
$$
\\

\noindent\textbf{Proof of Theorem 1}.  Let $G$ be a 1-tough graph. If $c\ge2\delta+4$ then we are done. Hence, we can assume that
$$
c\le2\delta+3.                                 \eqno{(1)}
$$

Let $C$ be a longest cycle in $G$ and $P=x_1\overrightarrow{P}x_2$ a longest path in $G\backslash C$. Put $|P|=|V(P)|-1=\overline{p}$. If $|V(P)|= 0$ then $C$ is a Hamilton cycle and we are done. Let $|V(P)|\ge1$, that is $\overline{p}\ge 0$.  Put $X=N_C(x_1)\cup N_C(x_2)$ and let $\xi_1,...,\xi_s$ be the elements of $X$ occuring on $C$ in a consecutive order. Put
$$
I_i=\xi_i\overrightarrow{C}\xi_{i+1},       \     I_i^\ast=\xi_i^+\overrightarrow{C}\xi_{i+1}^-  \  \ (i=1,...,s),
$$
where $\xi_{s+1}=\xi_1.$  Since $G$ is a 1-tough graph, we have $\delta\ge2$. \\

\textbf{Case 1}. $\overline{p}\le\delta-2$.

It follows that $s\ge|N_C(x_i)|\ge\delta-\overline{p}\ge2$ $(i=1,2)$. Assume first that $N_C(x_1)\not=N_C(x_2)$, implying that $\overline{p}\ge1$. If $\overline{p}\ge2$ then by Lemma 1, $c\ge4\delta-2\overline{p}\ge2\delta+4$, contradicting (1). Hence $\overline{p}=1$, which yields $\delta\ge\overline{p}+2=3$. By Lemma 1, $c\ge3\delta\ge9$. If $\delta\ge4$ then $c\ge3\delta\ge2\delta+4$, contradicting (1). Let $\delta=3$. Next, 
we can suppose that $c=9$, since otherwise $c\ge10=3\delta+1=2\delta+4$, contradicting (1).  Further, we can suppose that $s\ge3$, since  $N_C(x_1)=N_C(x_2)$ when $s=2$, contradicting the hypothesis. Finally, we can suppose that $s=3$, since clearly $c\ge10$ when $s\ge 4$, a contradiction. Thus, $|I_1|=|I_2|=|I_3|=3$ and it is not hard to see that $G\backslash \{\xi_1,\xi_2,\xi_3\}$ has at least four components, contradicting $\tau\ge1$.

Now assume that $N_C(x_1)=N_C(x_2)$. Since $C$ is extreme, we have 
$$
|I_i|\ge|\xi_ix_1\overrightarrow{P}x_2\xi_{i+1}|\ge\overline{p}+2 \ \ (i=1,...,s).
$$

\textbf{Case 1.1}. $s\ge\delta-\overline{p}+1$.

Clearly,
$$
c=\sum_{i=1}^s|I_i|\ge s(\overline{p}+2)
$$
$$
\ge(\delta-\overline{p}+1)(\overline{p}+2)=(\delta-\overline{p}-2)\overline{p}+2\delta+\overline{p}+2.  \eqno{(2)}
$$
If $\overline{p}\ge2$ then by (2), $c\ge2\delta+4$, contradicting (1). Let $\overline{p}\le1$.\\

\textbf{Case 1.1.1}. $\overline{p}=0$.

If $\Upsilon(I_1,...,I_s)=\emptyset$ then $G\backslash \{\xi_1,...,\xi_s\}$ has at least $s+1$ components, contradicting the fact that $\tau\ge1$. Otherwise   $\Upsilon(I_a,I_b)\not=\emptyset$ for some distinct $a,b\in \{1,...,s\}$.  Let $L\in \Upsilon(I_a,I_b)$. By Lemma 2(a1), 
$$
|I_a|+|I_b|\ge 2\overline{p}+2|L|+4\ge6.
$$
Recalling also that $s\ge\delta-\overline{p}+1=\delta+1$, we get
$$
c=\sum_{i=1}^s|I_i|\ge|I_a|+|I_b|+2(s-2)=2s+2\ge2\delta+4,
$$ 
contradicting (1).\\

\textbf{Case 1.1.2}. $\overline{p}=1$.

By (2), $c\ge3\delta$. We can suppose that $\delta\le3$, since $c\ge3\delta\ge2\delta+4$ when $\delta\ge4$, contradicting (1). On the other hand, by the hypothesis, $\delta\ge\overline{p}+2=3$, implying that $\delta=3$. By the hypothesis, $s\ge\delta-\overline{p}+1=3$. Next, we can suppose that $s=3$, since $c\ge s(\overline{p}+2)\ge12=2\delta+6$ when $s\ge4$, contradicting (1).  Further, if $\Upsilon(I_1,I_2,I_3)=\emptyset$ then $G\backslash \{\xi_1,\xi_2,\xi_3\}$ has at least four components, contradicting $\tau\ge1$. Otherwise $\Upsilon(I_a,I_b)\not=\emptyset$ for some distinct $a,b\in \{1,2,3\}$, say $a=1$ and $b=2$.  Let $L\in \Upsilon(I_1,I_2)$. By Lemma 2(a1), 
$$
|I_1|+|I_2|\ge 2\overline{p}+2|L|+4=8,
$$
which yields $c\ge|I_1|+|I_2|+|I_3|\ge11=2\delta+5$, contradicting (1).\\

\textbf{Case 1.2}. $s=\delta-\overline{p}$. 

It follows that $x_1x_2\in E$. Then $x_1x_2\overleftarrow{P}x_1^+$ is another longest path in $G\backslash C$. We can suppose that $N_C(x_1)=N_C(x_1^+)$, since otherwise we can argue as in Case 1. By the same reason,
$$
N_C(x_1)=N_C(x_1^+)=N_C(x_1^{+2})=...=N_C(x_2).
$$
Since $C$ is extreme, we have $|I_i|\ge|\xi_ix_1\overrightarrow{P}x_2\xi_{i+1}|=\overline{p}+2$ $(i=1,...,s)$. 
If $\Upsilon(I_1,...,I_s)=\emptyset$ then $G\backslash \{\xi_1,...,\xi_s\}$ has at least $s+1$ components, contradicting  $\tau\ge1$. Otherwise   $\Upsilon(I_a,I_b)\not=\emptyset$ for some distinct $a,b\in \{1,...,s\}$.  Let $L\in \Upsilon(I_a,I_b)$ with $L=z_1\overrightarrow{L}z_2$, where $z_1\in V(I_a^\ast)$ and $z_2\in V(I_b^\ast)$. By Lemma 2(a1), $|I_a|+|I_b|\ge2\overline{p}+6$. Hence
$$
c=\sum_{i=1}^s|I_i|\ge|I_a|+|I_b|+(s-2)(\overline{p}+2)\ge s(\overline{p}+2)+2
$$ 
$$
=(\delta-\overline{p})(\overline{p}+2)+2=2\delta+2+\overline{p}(\delta-\overline{p}-2).      \eqno{(3)}
$$

\textbf{Claim 1}. $(a1)$ $2\overline{p}+6\le |I_a|+|I_b|\le2\overline{p}+7$ and $|I_i|\le\overline{p}+5$  $(i=1,...,s)$.\

$(a2)$ If $|I_a|+|I_b|=2\overline{p}+7$ then  $|I_i|=\overline{p}+2$ for each $i\in \{1,...,s\}\backslash \{a,b\}$. 

$(a3)$ If $|I_a|+|I_b|=2\overline{p}+6$ then  $|I_f|\le\overline{p}+3$ for some $f\in \{1,...,s\}\backslash \{a,b\}$ and $|I_i|=\overline{p}+2$ for each $i\in \{1,...,s\}\backslash \{a,b,f\}$.  

$(a4)$ If $|I_f|=\overline{p}+5$ for some $f\in\{a,b\}$ then $|I_i|=\overline{p}+2$ for each $i\in \{1,...,s\}\backslash \{f\}$.

$(a5)$ For each distinct $f,g,h\in\{1,...,s\}$, $|I_f|+|I_g|+|I_h|\le3\overline{p}+9$.

$(a6)$ $\Upsilon(I_1,...,I_s)\subseteq E$.

\textbf{Proof}. If $|I_f|\ge\overline{p}+6$ for some $f\in \{1,...,s\}$ then 
$$
c=\sum_{i=1}^s|I_i|\ge|I_f|+(s-1)(\overline{p}+2)\ge s(\overline{p}+2)+4
$$ 
$$
=2\delta+4+\overline{p}(\delta-\overline{p}-2)\ge 2\delta+4, 
$$
contradicting (1). Next, if  $|I_a|+|I_b|\ge2\overline{p}+8$  then 
$$
c\ge|I_a|+|I_b|+(s-2)(\overline{p}+2)\ge s(\overline{p}+2)+4\ge2\delta+4,
$$ 
again contradicting (1). Hence $(a1)$ holds. Statements $(a2)-(a4)$ can be proved by a similar way. To prove $(a5)$, assume the contrary, that is $|I_f|+|I_g|+|I_h|\ge3\overline{p}+10$ for some distinct $f,g,h\in\{1,...,s\}$. Then 
$$
c=\sum_{i=1}^s|I_i|\ge|I_f|+|I_g|+|I_h|+(s-3)(\overline{p}+2)
$$
$$
\ge 3(\overline{p}+2)+4+(s-3)(\overline{p}+2)=2\delta+4+\overline{p}(s-2)\ge2\delta+4,
$$
contradicting (1). Statement $(a6)$ follows from Lemma 2(a1) and Claim 1(a1). Claim 1 is proved.\\

\textbf{Claim 2}. $\overline{p}+3\le d_1\le\overline{p}+4$ and $\overline{p}+3\le d_2\le\overline{p}+4$, where
$$
d_1=|\xi_a\overrightarrow{C}z_1|+|\xi_b\overrightarrow{C}z_2|, \ \ d_2=|z_1\overrightarrow{C}\xi_{a+1}|+|z_2\overrightarrow{C}\xi_{b+1}|.
$$

\textbf{Proof}. Put 
$$
Q=\xi_ax_1\overrightarrow{P}x_2\xi_b\overleftarrow{C}z_1z_2\overrightarrow{C}\xi_a.
$$
Clearly, $|Q|=|C|-d_1+\overline{p}+3$. Since $C$ is extreme, we have $|C|\ge|Q|$, implying that $d_1\ge \overline{p}+3$. By a symmetric argument, $d_2\ge \overline{p}+3$. By Claim 1(a1), $|I_a|+|I_b|=d_1+d_2\le2\overline{p}+7$. 
If $d_1\ge \overline{p}+5$ then 
$2\overline{p}+7\ge d_1+d_2\ge\overline{p}+5+d_2$, implying that 
$d_2\le \overline{p}+2$, a contradiction. Hence, 
$d_1\le \overline{p}+4$. By a symmetric argument, 
$d_2\le \overline{p}+4$. Claim 2 is proved.\\

\textbf{Claim 3}. If $v_1\in V(\xi_a^+\overrightarrow{C}z_1^-)$ and $v_2\in V(z_1^+\overrightarrow{C}\xi_{a+1}^-)$ then $v_1v_2\not\in E$.

\textbf{Proof}. Assume the contrary, that is $v_1v_2\in E$. Put

$$
Q=\xi_a\overrightarrow{C}v_1v_2\overleftarrow{C}z_1z_2\overleftarrow{C}\xi_{a+1}x_1\overrightarrow{P}x_2\xi_{b+1}\overrightarrow{C}\xi_a,
$$
$$
|\xi_a\overrightarrow{C}v_1|=d_1, \ \  |v_1\overrightarrow{C}z_1|=d_2,  \ \ |z_1\overrightarrow{C}v_2|=d_3, 
$$
$$
|v_2\overrightarrow{C}\xi_{a+1}|=d_4, \ \  |\xi_b\overrightarrow{C}z_2|=d_5, \ \  |z_2\overrightarrow{C}\xi_{b+1}|=d_6.
$$
Clearly, $|Q|=|C|-d_2-d_4-d_6+\overline{p}+4$. Since $C$ is extreme, we have $|Q|\le |C|$, implying that $d_2+d_4+d_6\ge\overline{p}+4$. By a symmetric argument, $d_1+d_3+d_5\ge\overline{p}+4$. By summing, we get
$$
\sum_{i=1}^6d_i=|I_a|+|I_b|\ge2\overline{p}+8,
$$
contradicting Claim 1(a1). Thus, $v_1v_2\not\in E$. Claim 3 is proved.\\

\textbf{Claim 4}. Let $\xi_f,\xi_g,\xi_h$ occur on $\overrightarrow{C}$ in a consecutive order for some $f,g,h\in \{1,...,s\}$ and $w_1w_2\in E$ for some $w_1\in V(I_f^\ast)$ and $w_2\in V(I_g^\ast)$. If $N(w_3)\cap \{\xi_{f+1},\xi_g\}\not=\emptyset$ for some $w_3\in V(I_h^\ast)$ then
$$
|w_1\overrightarrow{C}\xi_{f+1}|+|\xi_{g}\overrightarrow{C}w_2|+|\xi_{h}\overrightarrow{C}w_3|\ge\overline{p}+4.
$$ 
Further, if $N(w_4)\cap \{\xi_{f+1},\xi_g\}\not=\emptyset$ for some $w_4\in V(I_{h-1}^\ast)$ then 
$$
|w_1\overrightarrow{C}\xi_{f+1}|+|\xi_{g}\overrightarrow{C}w_2|+|w_4\overrightarrow{C}\xi_h|\ge\overline{p}+4.
$$
\textbf{Proof}. Assume first that $w_3\xi_{f+1}\in E$. Put
$$
Q=\xi_f\overrightarrow{C}w_1w_2\overrightarrow{C}\xi_hx_1\overrightarrow{P}x_2\xi_g\overleftarrow{C}\xi_{f+1}w_3\xi_f.
$$ 
Clearly,
$$
|Q|=|C|-|w_1\overrightarrow{C}\xi_{f+1}|-|\xi_g\overrightarrow{C}w_2|-|\xi_h\overrightarrow{C}w_3|+\overline{p}+4.
$$
Since $|Q|\le |C|$, the desired result holds immediately. If $w_4\xi_{f+1}\in E$ then we can use the following cycle
$$
Q^\prime=\xi_f\overrightarrow{C}w_1w_2\overrightarrow{C}w_4\xi_{f+1}\overrightarrow{C}\xi_gx_2\overleftarrow{P}x_1\xi_h\overrightarrow{C}\xi_f
$$
instead of $Q$. By a symmetric argument, the desired result holds when either $w_3\xi_g\in E$ or $w_4\xi_g\in E$. Claim 4 is proved.\\

\textbf{Claim 5}. Every two intermediate independent edges $e_1,e_2$ in $\Upsilon(I_1,...,I_s)$ have a crossing with $e_1,e_2\in \Upsilon(I_f,I_g,I_h)$ for some distinct $f,g,h\in\{1,...,s\}$.  

\textbf{Proof}. Let $e_1=w_1w_2$ and $e_2=w_3w_4$. We distinguish three different cases. First, if $e_1,e_2\in \Upsilon(I_f,I_g)$ for some distinct $f,g$, then by Lemma 2(a3), $|I_f|+|I_g|\ge2\overline{p}+8$, contradicting Claim 1(a1). Next, if  $e_1\in \Upsilon(I_f,I_g)$ and $e_2\in \Upsilon(I_h,I_r)$ for some distinct $f,g,h,r$, then by Lemma 2(a1), $|I_f|+|I_g|\ge2\overline{p}+6$ and $|I_h|+|I_r|\ge2\overline{p}+6$, implying that
$$
c\ge|I_f|+|I_g|+|I_h|+|I_r|+(s-4)(\overline{p}+2)=4\overline{p}+12+(s-4)(\overline{p}+2)
$$
$$
=s(\overline{p}+2)+4=2\delta+4+\overline{p}(\delta-\overline{p}-2)\ge2\delta+4,
$$
which again contradicts (1). Finally, let $e_1\in \Upsilon(I_f,I_g)$ and $e_2\in \Upsilon(I_f,I_h)$ for some distinct $f,g,h$. Assume w.l.o.g. that $\xi_f,\xi_g,\xi_h$ occur on $\overrightarrow{C}$ in a consecutive order and $w_1,w_3\in V(I_f^\ast)$, $w_2\in V(I_g^\ast)$, $w_4\in V(I_h^\ast)$. We can assume also that $w_3$ and $w_1$ occur on $I_f$ in a consecutive order, since otherwise $e_1$ and $e_2$ have a crossing and we are done. Put
$$
Q=\xi_f\overrightarrow{C}w_3w_4\overleftarrow{C}w_2w_1\overrightarrow{C}\xi_gx_2\overleftarrow{P}x_1\xi_{h+1}\overrightarrow{C}\xi_f,
$$
$$
|\xi_f\overrightarrow{C}w_3|=d_1, \  |w_3\overrightarrow{C}w_1|=d_2,  \  |w_1\overrightarrow{C}\xi_{f+1}|=d_3,
$$
$$
|\xi_g\overrightarrow{C}w_2|=d_4, \  |w_2\overrightarrow{C}\xi_{g+1}|=d_5,  \  |\xi_h\overrightarrow{C}w_4|=d_6
,  \  |w_4\overrightarrow{C}\xi_{h+1}|=d_7.
$$
Clearly, $|Q|=|C|-d_2-d_4-d_7+\overline{p}+4$. Since $C$ is extreme, we have $|Q|\le |C|$, implying that $d_2+d_4+d_7\ge\overline{p}+4$. On the other hand, by Lemma 2, $d_3+d_5\ge \overline{p}+3$ and $d_1+d_6\ge \overline{p}+3$. By summing, we get $\sum_{i=1}^7d_i=|I_f|+|I_g|+|I_h|\ge3\overline{p}+10$. Then
$$
|C|\ge|I_f|+|I_g|+|I_h|+(s-3)(\overline{p}+2)=s(\overline{p}+2)+4\ge2\delta+4,
$$ 
contradicting (1). Claim 5 is proved.\\

\textbf{Claim 6}. If $\mu(\Upsilon)=1$ then $s\le3$ and either $\xi_a^+\xi_{b+1}^-\in E$ with $\xi_a=\xi_{b+1}$ or $\xi_{a+1}^-\xi_b^+\in E$ with $\xi_{a+1}=\xi_b$. If $\mu(\Upsilon)=1$ and $s=3$ then $|I_1|=|I_2|=|I_3|=\overline{p}+3$.

\textbf{Proof}. Since $\mu(\Upsilon)=1$, either one of the vertices $z_1, z_2$, say $z_1$, is a common vertex for all edges in $\Upsilon(I_1,...,I_s)$ or $z_1z_3,z_2z_3\in \Upsilon(I_1,...,I_s)$ for some $z_3\in V(I_f^\ast)$ and $f\in\{1,...,s\}\backslash \{a,b\}$.\\

\textbf{Case a1}. $z_1$ is a common vertex for all edges in $\Upsilon(I_1,...,I_s)$.

If $z_1\not\in \{\xi_a^+,\xi_{a+1}^-\}$ then by Claim 3, $G\backslash \{\xi_1,...,\xi_s,z_1\}$ has at least $s+2$ components, contradicting $\tau\ge1$. Let $z_1\in \{\xi_a^+,\xi_{a+1}^-\}$, say $z_1=\xi_a^+$.\\

\textbf{Case a1.1}. $z_1\xi_{b+1}^-\not\in E$.

It follows that $z_2\not=\xi_{b+1}^-$. By Claim 2, $|\xi_b\overrightarrow{C}z_2|\ge\overline{p}+2$.\\

\textbf{Case a1.1.1}. $z_1\xi_{b+1}^{-2}\not\in E$.

It follows that $|I_b|\ge \overline{p}+5$. By Claim 1(a1), $|I_a|=\overline{p}+2$. Moreover, we have $|I_b|= \overline{p}+5$, $|\xi_b\overrightarrow{C}z_2|=\overline{p}+2$, $z_2=\xi_{b+1}^{-3}$ and $N(z_1)\cap V(I_b^\ast)=\{z_2\}$. By Claim 1(a4), $|I_i|=\overline{p}+2$ for each $i\in\{1,...,s\}\backslash \{b\}$. Next, by Lemma 2(a1), $\Upsilon(I_a,I_i)=\emptyset$ for each $i\in\{1,...,s\}\backslash \{a,b\}$. Thus, if $z_1y\in \Upsilon(I_1,...,I_s)$ then $y=z_2$, implying that $\Upsilon(I_1,...,I_s)=\{z_1z_2\}$. Besides,  since $|\xi_b\overrightarrow{C}z_2|=\overline{p}+2\ge2$, we have $z_2\not\in\{\xi_b^+,\xi_{b+1}^-\}$. Therefore, by Claim 3,  $G\backslash \{\xi_1,...,\xi_s,z_2\}$ has at least $s+2$ components, contradicting $\tau\ge1$. \\

\textbf{Case a1.1.2}. $z_1\xi_{b+1}^{-2}\in E$.

It follows that $|I_b|\ge\overline{p}+4$. Assume first that $|I_b|=\overline{p}+5$. If $z_1\xi_{b+1}^{-3}\not\in E$ then clearly $z_2=\xi_{b+1}^{-2}$ and we can argue as in Case a1.1.1. Otherwise the following cycle
$$
\xi_ax_1\overrightarrow{P}x_2\xi_{a+1}\overrightarrow{C}\xi_{b+1}^{-3}z_1\xi_{b+1}^{-2}\overrightarrow{C}\xi_a
$$
is longer than $C$, a contradiction. 

Now assume that $|I_b|=\overline{p}+4$, that is $|\xi_b\overrightarrow{C}\xi_{b+1}^{-2}|=\overline{p}+2$. If $z_1y\in E$ for some $y\in V(\xi_b\overrightarrow{C}\xi_{b+1}^{-3})$ then by Claim 2, $|\xi_b\overrightarrow{C}y|\ge \overline{p}+2$, implying that  $|I_b|\ge\overline{p}+5$, a contradiction. Hence, if $z_1y\in \Upsilon(I_a,I_b)$ then clearly $y=\xi_{b+1}^{-2}$. In particular, we have $z_2=\xi_{b+1}^{-2}$. Further, if $z_1y\in \Upsilon(I_a,I_f)$ for some $f\in\{1,...,s\}\backslash\{b\}$, then by Lemma 2(a1), $|I_a|+|I_f|\ge2\overline{p}+6$, that is $|I_a|+|I_b|+|I_f|\ge3\overline{p}+10$, contradicting Claim 1(a5). Thus, $z_2$ is a common vertex for all edges in $\Upsilon(I_1,...,I_s)$. By Claim 3, $G\backslash \{\xi_1,...,\xi_s,z_2\}$ has at least $s+2$ components, contradicting $\tau\ge1$.\\

\textbf{Case a1.2}. $\xi_a^+\xi_{b+1}^-\in E$.

By Claim 2, $|\xi_a^+\overrightarrow{C}\xi_{a+1}|\ge\overline{p}+2$ and $|\xi_b\overrightarrow{C}\xi_{b+1}^-|\ge\overline{p}+2$. If $|\xi_a^+\overrightarrow{C}\xi_{a+1}|\ge\overline{p}+3$ and $|\xi_b\overrightarrow{C}\xi_{b+1}^-|\ge\overline{p}+3$ then $|I_a|+|I_b|\ge2\overline{p}+8$, contradicting Claim 1(a1). Hence, we can assume w.l.o.g. that $|\xi_b\overrightarrow{C}\xi_{b+1}^-|=\overline{p}+2$, that is $|I_b|=\overline{p}+3$ and $|I_a|\ge\overline{p}+3$. Further, we have $\xi_b^+\xi_a,\xi_b^+\xi_{b+1}\not\in E$ (by Claim 4) and $\xi_b^+\xi_a^+\not\in E$ (by Claim 2).\\

\textbf{Case a1.2.1}. $N(\xi_b^+)\not\subseteq V(C)$.

Let $Q=\xi_b^+\overrightarrow{Q}v$ be a longest path in $G$ with $V(Q)\cap V(C)=\{\xi_b^+\}$.  Since $C$ is extreme, we have $V(Q)\cap V(P)=\emptyset$. Next, since $P$ is a longest path in $G\backslash C$, we have $|Q|\le\overline{p}+1$. Further, recalling that $\xi_b^+\xi_a,\xi_b^+\xi_{b+1},\xi_b^+\xi_a^+\not\in E$ (see Case a1.2), we conclude that $v\xi_a, v\xi_{b+1}, v\xi_a^+\not\in E$, as well. If $vy\not\in E$ for each $y\in (\xi_b^{+2}\overrightarrow{C}\xi_{b+1}^-)$ then clearly
$$
N(v)\subseteq(V(Q)\cup\{\xi_1,...,\xi_s\})\backslash\{\xi_a,\xi_{b+1}\xi_a^+\},
$$ 
that is $d(v)\le|Q|+s-2\le\overline{p}+s-1=\delta-1$, a contradiction. Now let $vy\in E$ for some $y\in V(\xi_b^{+2}\overrightarrow{C}\xi_{b+1}^-)$. Assume that $y$ is chosen so as to minimize $|\xi_b^+\overrightarrow{C}y|$. Since $C$ is extreme, we have $|\xi_b^+\overrightarrow{C}y|\ge|Q|+1$. Further, since 
$$
|N(v)\cap V(y\overrightarrow{C}\xi_{b+1}^-)|\ge\delta-(s-2)-|Q|,
$$
we have 
$$
|\xi_b^+\overrightarrow{C}\xi_{b+1}^-|\ge|Q|+1+2(\delta-s+1-|Q|)
$$
$$
=2\delta-|Q|-2s+3\ge2\delta-\overline{p}-2s+2=\overline{p}+2.
$$
But then $|I_b|\ge\overline{p}+4$, a contradiction.\\

\textbf{Case a1.2.2}. $N(\xi_b^+)\subseteq V(C)$.

Since $\mu(\Upsilon)=1$ and $\xi_b^+\xi_a^+\not\in E$, we have 
$$
N(\xi_b^+)\subseteq V(\xi_b^{+2}\overrightarrow{C}\xi_{b+1}^-)\cup \{\xi_1,...,\xi_s\}\backslash \{\xi_a,\xi_{b+1}\}.
$$
If $\xi_a\not=\xi_{b+1}$ then $d(\xi_b^+)\le\overline{p}+s-1=\delta-1$, a contradiction. Hence $\xi_a=\xi_{b+1}$. \\

\textbf{Case a1.2.2.1}. $|I_f|=\overline{p}+2$ for some $f\in\{1,...,s\}\backslash \{a,b\}$. 

If $N(\xi_f^+)\subseteq V(C)$ then as above, 
$$
d(\xi_f^+)\le s-1+|\xi_f^+\overrightarrow{C}\xi_{f+1}^-|=\overline{p}+s-1=\delta-1,
$$
a contradiction. If $N(\xi_f^+)\not\subseteq V(C)$ then we can argue as in Case a1.2.1. \\

\textbf{Case a1.2.2.2}. $|I_i|\ge\overline{p}+3$ for each $i\in\{1,...,s\}\backslash \{a,b\}$. 

If $s\ge4$ then
$$
|C|=\sum_{i=1}^s|I_i|\ge s(\overline{p}+3)=(\delta-\overline{p})(\overline{p}+3)
$$
$$
=2\delta+2\overline{p}+4+(\delta-\overline{p}-4)(\overline{p}+1)\ge2\delta+4,
$$
contradicting (1). Hence, $s\le3$. Moreover, if $s=3$ then by Claim 1(a5), $|I_1|=|I_2|=|I_3|=\overline{p}+3$.\\

\textbf{Case a2}. $z_1z_3,z_2z_3\in \Upsilon(I_1,...,I_s)$, where $z_3\in V(I_f^\ast)$ and $f\in\{1,...,s\}\backslash \{a,b\}$.

Assume w.l.o.g. that $\xi_a,\xi_b,\xi_f$ occur on $\overrightarrow{C}$ in a consecutive order. Put
$$
|\xi_a\overrightarrow{C}z_1|=d_1, \ |z_1\overrightarrow{C}\xi_{a+1}|=d_2, \ |\xi_b\overrightarrow{C}z_2|=d_3,
$$
$$
|z_2\overrightarrow{C}\xi_{b+1}|=d_4, \ |\xi_f\overrightarrow{C}z_3|=d_5, \ |z_3\overrightarrow{C}\xi_{f+1}|=d_6.
$$
By Claim 2, 
$$
d_1+d_3\ge\overline{p}+3, \  d_1+d_5\ge\overline{p}+3, \  d_2+d_4\ge\overline{p}+3, 
$$
$$
d_2+d_6\ge\overline{p}+3, \  d_3+d_5\ge\overline{p}+3, \  d_4+d_6\ge\overline{p}+3.
$$
By summing, we get
$$
2\sum_{i=1}^6d_i=2(|I_a|+|I_b|+|I_f|)\ge6(\overline{p}+3).
$$
On the other hand, by Claim 1(a5), $|I_a|+|I_b|+|I_f|\le3(\overline{p}+3)$, implying that $d_1=d_2=...=d_6=(\overline{p}+3)/2$ and $\overline{p}$ is odd. Hence $d_i\ge2$ and using Claim 3, we can state that $G\backslash \{\xi_1,...,\xi_s,z_1,z_2\}$ has at least $s+3$ components, contradicting $\tau\ge1$. Claim 6 is proved.\\

\textbf{Claim 7}. Either $\mu(\Upsilon)=1$ or $\mu(\Upsilon)=3$.

\textbf{Proof}. The proof is by contradiction. If $\mu(\Upsilon)=0$ then $G\backslash \{\xi_1,...,\xi_s\}$ has at least $s+1$ components, contradicting $\tau\ge1$. Let $\mu(\Upsilon)\ge 1$.\\

\textbf{Case a1}. $\mu=2$.

By Claim 5, $\Upsilon(I_1,...,I_s)$ consists of two crossing intermediate independent edges $w_1w_2\in \Upsilon(I_f,I_g)$ and $w_3w_4\in \Upsilon(I_f,I_h)$ for some distinct $f,g,h$. Assume that both $\xi_f,\xi_g,\xi_h$ and $w_1,w_3,w_2,w_4$  occur on $\overrightarrow{C}$ in a consecutive order. Put
$$
Q=\xi_f\overrightarrow{C}w_1w_2\overrightarrow{C}w_4w_3\overrightarrow{C}\xi_gx_2\overleftarrow{P}x_1\xi_{h+1}\overrightarrow{C}\xi_f,
$$
$$
|\xi_f\overrightarrow{C}w_1|=d_1, \ |w_1\overrightarrow{C}w_3|=d_2,  \  |w_3\overrightarrow{C}\xi_{f+1}|=d_3,
$$
$$
|\xi_g\overrightarrow{C}w_2|=d_4, \  |w_2\overrightarrow{C}\xi_{g+1}|=d_5, \  |\xi_h\overrightarrow{C}w_4|=d_6, \  |w_4\overrightarrow{C}\xi_{h+1}|=d_7.
$$
Clearly, $|Q|=|C|-d_2-d_4-d_7+\overline{p}+4$. Since $|Q|\le |C|$, we have $d_2+d_4+d_7\ge\overline{p}+4$. If $d_3+d_6\ge\overline{p}+3$ and $d_1+d_5\ge\overline{p}+3$ then $\sum_{i=1}^7d_i=|I_f|+|I_g|+|I_h|\ge3\overline{p}+10$, contradicting Claim 1(a5). Otherwise, either $d_3+d_6\le\overline{p}+2$ or $d_1+d_5\le\overline{p}+2$, say  $d_3+d_6\le\overline{p}+2$. Further, if either $d_7=1$ or $\xi_{h+1}^-w_3\in E$ then by Claim 2, $d_3\ge \overline{p}+2$, that is $d_3+d_6\ge\overline{p}+3$, a contradiction. Hence, $d_7\ge2$ and $\xi_{h+1}^-w_3\not\in E$. By Claim 4, $\xi_{h+1}^-\xi_{f+1}, \xi_{h+1}^-\xi_{h}\not\in E$. If $|I_h|\ge \overline{p}+4$ then taking account that $|I_f|+|I_g|\ge2\overline{p}+6$ (by Claim 1(a1)), we get $|I_f|+|I_g|+|I_h|\ge3\overline{p}+10$, contradicting Claim 1(a5). Hence, $|I_h|\le\overline{p}+3$. By a symmetric argument,  $|I_g|\le\overline{p}+3$.\\

\textbf{Case a1.1}. $N(\xi_{h+1}^-)\subseteq V(C)$.

If $\xi_{h+1}^-w_2\not\in E$ then recalling that $\mu(\Upsilon)=2$, we get
$$
N(\xi_{h+1}^-)\subseteq V(w_4\overrightarrow{C}\xi_{h+1}^{-2})\cup\{\xi_1,...,\xi_s\}\backslash\{\xi_{f+1},\xi_h\},
$$
implying that $|N(\xi_{h+1}^-)|\le\overline{p}+s-1=\delta-1$, a contradiction. Now let $\xi_{h+1}^-w_2\in E$. By Claim 1(a1 and a5), $|I_f|=|I_g|=|I_h|=\overline{p}+3$. Moreover, by Claim 2, $d_5=\overline{p}+2$ and $d_4=1$. Then by, by the same reason, $d_1=\overline{p}+2$, implying that $|I_a|\ge\overline{p}+4$, a contradiction.\\

\textbf{Case a1.2}. $N(\xi_{h+1}^-)\not\subseteq V(C)$.

We can argue as in the proof of Claim 6 (Case a1.2.1).\\

\textbf{Case a2}. $\mu(\Upsilon)\ge4$.

By Claim 5, there are at least four pairwise crossing intermediate independent edges in $\Upsilon(I_1,...,I_s)$, which is impossible. Claim 7 is proved.\\

\textbf{Claim 8}. If $\mu(\Upsilon)=1$ then either $n\equiv1(mod\ 3)$ with $c\ge2\delta+2$ or $n\equiv1(mod\ 4)$ with $c\ge2\delta+3$ or $n\equiv2(mod\ 3)$ with $c\ge2\delta+3$.

  \textbf{Proof}. By Claim 6, $s\le3$ and either $\xi_a^+\xi_{b+1}^-\in E$ or $\xi_{a+1}^-\xi_b^+\in E$, say $\xi_{a+1}^-\xi_b^+\in E$. \\

\textbf{Case a1}. $s=2$.

It follows that $\delta=\overline{p}+s=\overline{p}+2$. Let $a=1$ and $b=2$. By Claim 2, $|\xi_1\overrightarrow{C}\xi_2^-|\ge \overline{p}+2$ and $|\xi_2^+\overrightarrow{C}\xi_1|\ge \overline{p}+2$, implying that $|I_i|\ge \overline{p}+3$ $(i=1,2)$.\\

 \textbf{Case a1.1}. $|I_1|=\overline{p}+4$ and $|I_2|=\overline{p}+3$.

If $V(G)=V(C\cup P)$ then $n=3\overline{p}+8=3\delta+2\equiv2(mod\ 3)$ with $c=2\overline{p}+7=2\delta+3$, and we are done. Otherwise $N(v_1)\not\subseteq V(C\cup P)$ for some $v_1\in V(C\cup P)$. Observing that $x_1x_2\in E$ and recalling that $P$ is a longest path in $V(G\backslash C)$, we conclude that  $v_1\not\in V(P)$. Choose a longest path $Q=v_1\overrightarrow{Q}v_2$ with $V(Q)\cap V(C)=\{v_1\}$. Clearly, $1\le|Q|\le\overline{p}+1=\delta-1$ and $N(v_2)\subseteq V(C\cup Q)$.\\

\textbf{Case a1.1.1}. $v_1\in V(\xi_2^{+2}\overrightarrow{C}\xi_1^-)$.

By Claim 1(a6), $N(v_2)\cap V(I_1^\ast)=\emptyset$, that is $N(v_2)\subseteq V(I_1)\cup V(Q)$. Assume that $v_1$ is chosen so as to minimize $|v_1\overrightarrow{C}\xi_1|$, implying that $N(v_2)\cap V(v_1\overrightarrow{C}\xi_1^-)=\emptyset$. Clearly, $|v_1\overrightarrow{C}\xi_1|\le \overline{p}+1$. Then by Claim 4, $v_1\xi_2\not\in E$ and therefore,  $v_2\xi_2\not\in E$, as well. \\

\textbf{Case a1.1.1.1}. $v_2\xi_1\in E$.

It follows that $N(v_2)\subseteq V(Q)\cup V(\xi_2^+\overrightarrow{C}v_1^-)\cup \{\xi_1\}$. Since $C$ is extreme and $v_2\xi_1\in E$, we have $|v_1\overrightarrow{C}\xi_1|\ge |Q|+1$. If $N(v_2)\subseteq V(Q)\cup \{\xi_1\}$ then clearly $|Q|\ge\delta-1=\overline{p}+1$ and therefore, $|v_1\overrightarrow{C}\xi_1|\ge \overline{p}+2$. But then $|I_2|\ge\overline{p}+4$, a contradiction. Hence, $N(v_2)\not\subseteq V(Q)\cup \{\xi_1\}$, that is $v_2y\in E$ for some $y\in V(\xi_2^+\overrightarrow{C}v_1^-)$. Assume that $y$ is chosen so as to minimize $|y\overrightarrow{C}v_1|$. Observing that $|y\overrightarrow{C}v_1|\ge |Q|+1$ and $\delta=|\xi_2^+\overrightarrow{C}\xi_1|\ge4$, we get
$$
|\xi_2^+\overrightarrow{C}\xi_1|\ge2(|Q|+1)+2(\delta-|Q|-2)=2\delta-2\ge\delta+2=\overline{p}+4,
$$
a contradiction. \\

\textbf{Case a1.1.1.2}. $v_2\xi_1\not\in E$.

It follows that $N(v_2)\subseteq V(Q)\cup V(\xi_2^+\overrightarrow{C}v_1^-)$. If $N(v_2)\subseteq V(Q)$ then  $|Q|\ge\delta=\overline{p}+2$, a contradiction. Otherwise $v_2y\in E$ for some  $y\in V(\xi_2^+\overrightarrow{C}v_1^-)$.  Assume that $y$ is chosen so as to minimize $|y\overrightarrow{C}v_1|$. Since $|y\overrightarrow{C}v_1|\ge |Q|+1$, we have
$$
|\xi_2^+\overrightarrow{C}v_1|\ge|Q|+1+2(\delta-|Q|-1)=2\delta-|Q|-1\ge\delta=\overline{p}+2.
$$
But then $|I_b|\ge4$, a contradiction.\\

\textbf{Case a1.1.2}. $v_1\in V(\xi_1^{+}\overrightarrow{C}\xi_2^{-3})$.

By Claim 1(a6), $N(v_2)\cap V(I_2^\ast)=\emptyset$, that is $N(v_2)\subseteq V(Q)\cup V(I_1)$. Assume that $v_1$ is chosen so as to minimize $|\xi_1\overrightarrow{C}v_1|$, implying that $N(v_2)\cap V(\xi_1^+\overrightarrow{C}v_1^-)=\emptyset$. Clearly, $|\xi_1\overrightarrow{C}v_1|\le \overline{p}+1$. Then by Claim 4, $v_1\xi_2\not\in E$ and therefore, $v_2\xi_2\not\in E$. \\

\textbf{Case a1.1.2.1}. $\xi_2^+\xi_2^{-2}\in E$.

By Claim 3, $v_1\xi_2^-\not \in E$, implying that $v_2\xi_2^-\not\in E$.\\

\textbf{Case a1.1.2.1.1}. $v_2\xi_1\in E$.

It follows that  $N(v_2)\subseteq V(Q)\cup V(v_1\overrightarrow{C}\xi_2^{-2})\cup \{\xi_1\}$. Since $C$ is extreme and $v_2\xi_1\in E$, we have $|\xi_1\overrightarrow{C}v_1|\ge |Q|+1$. If $N(v_2)\subseteq V(Q)\cup \{\xi_1\}$ then $|Q|\ge \delta-1=\overline{p}+1$ and therefore, $|\xi_1\overrightarrow{C}v_1|\ge \overline{p}+2$. But then $|I_1|\ge \overline{p}+5$, a contradiction. Hence, $N(v_2)\not\subseteq V(Q)\cup \{\xi_1\}$, that is $v_2y\in E$ for some $y\in V(v_1^+\overrightarrow{C}\xi_2^{-2})$. Assume that $y$ is chosen so as to minimize $|v_1\overrightarrow{C}y|$. Observing that $|v_1\overrightarrow{C}y|\ge|Q|+1$ and $\delta=|\xi_1\overrightarrow{C}\xi_2^{-2}|\ge4$, we get
$$
|\xi_1\overrightarrow{C}\xi_2^{-2}|\ge 2(|Q|+1)+2(\delta-|Q|-2)=2\delta-2\ge\delta+2=\overline{p}+4,
$$
a contradiction.\\

\textbf{Case a1.1.2.1.2}. $v_2\xi_1\not\in E$.

It follows that  $N(v_2)\subseteq V(Q)\cup V(v_1\overrightarrow{C}\xi_2^{-2})$.  If $N(v_2)\subseteq V(Q)$ then $|Q|\ge \delta=\overline{p}+2$, a contradiction. Otherwise  $v_2y\in E$ for some $y\in V(v_1^+\overrightarrow{C}\xi_2^{-2})$. By choosing $y$ so as to minimize $|v_1\overrightarrow{C}y|$, we get
$$
|v_1\overrightarrow{C}\xi_2^{-2}|\ge |Q|+1+2(\delta-|Q|-1)=2\delta-|Q|-1\ge\delta=\overline{p}+2.
$$
This yields $|I_a|\ge\overline{p}+5$, a contradiction.\\

\textbf{Case a1.1.2.2}. $\xi_2^+\xi_2^{-2}\not\in E$.

If $v_2\xi_1\in E$ then as in Case a1.1.2.1.1, $|\xi_1\overrightarrow{C}\xi_2^-|\ge \overline{p}+4$, contradicting  the fact that $|I_1|=\overline{p}+4$. Otherwise, as in Case a1.1.2.1.2, $|v_1\overrightarrow{C}\xi_2^-|\ge \overline{p}+2$. Since $|I_1|=\overline{p}+4$, we have $v_1=\xi_1^+$, $|Q|=\delta-1=\overline{p}+1$ and $v_3=\xi_2^-$. Moreover, we have $N(v_2)=(V(Q)\cup \{\xi_2^-\})\backslash \{v_2\}$. Further, let $v$ be an arbitrary vertex in $V(Q)\backslash \{v_1\}$. Put $Q^\prime=v_1\overrightarrow{Q}v^-v_2\overleftarrow{Q}v$. Since $Q^\prime$ is another longest path with $V(Q^\prime)\cap V(C)=\{v_1\}$, we can suppose that $N(v)=(V(Q)\cup \{\xi_2^-\})\backslash \{v\}$ for each $v\in V(Q)\backslash\{v_1\}$. Furthermore, if $\xi_1y\in E$ for some $y\in V(\xi_1^{+2}\overrightarrow{C}\xi_2^{-2})$ then
$$
\xi_1x_1\overrightarrow{P}x_2\xi_2\xi_2^+\xi_2^-v_2\overleftarrow{Q}v_1\overrightarrow{C}y\xi_1
$$
is longer than $C$, a contradiction. Hence, $\xi_1y\not\in E$ for each $y\in V(\xi_1^{+2}\overrightarrow{C}\xi_2^{-2})$. Analogously,  if $y\xi_2\in E$ for some $y\in V(\xi_1^{+}\overrightarrow{C}\xi_2^{-2})$ then
$$
\xi_1x_1\overrightarrow{P}x_2\xi_2y\overleftarrow{C}\xi_1^+\overrightarrow{Q}v_2\xi_2^-\xi_2^+\overrightarrow{C}\xi_1
$$
is longer than $C$, a contradiction. Hence, $y\xi_2\not\in E$ for each $y\in V(\xi_1^{+}\overrightarrow{C}\xi_2^{-2})$. But then $G\backslash \{\xi_1^+,\xi_2^-\}$ has at least three components, contradicting $\tau\ge1$.\\

 \textbf{Case a1.1.3}. $v_1=\xi_2^{-2}$.

By Claim 1(a6), $N(v_2)\subseteq V(I_1)$. If $v_2y\in E$ for some $y\in V(\xi_1^+\overrightarrow{C}v_1^{-})$ then we can argue as in Case a1.1.2. Hence, $N(v_2)\subseteq V(Q)\cup \{\xi_1,\xi_2\}$. If $v_2\xi_2\in E$ then
$$
\xi_1x_1\overrightarrow{P}x_2\xi_2v_2\overleftarrow{Q}v_1\xi_2^-\xi_2^+\overrightarrow{C}\xi_1
$$
is longer than $C$, a contradiction. Then clearly, $v_2\xi_1\in E$ and $N(v_2)\subseteq V(Q)\cup \{\xi_1\}$. Furthermore, we have $|Q|\ge \delta-1$, implying that $|\xi_1\overrightarrow{C}v_1|\ge|Q|+1\ge\delta$. Since $|\xi_1\overrightarrow{C}v_1|=\delta$, we have $|Q|=\delta-1=\overline{p}+1$ and $N(v_2)=(V(Q)\cup \{\xi_1\})\backslash\{v_2\}$. Moreover, as in Case 1.1.2.2, we have  $N(v)=(V(Q)\cup \{\xi_1\})\backslash\{v\}$ for each $v\in V(Q)\backslash \{v_1\}$. Now consider an arbitrary vertex $y\in V(\xi_1^+\overrightarrow{C}\xi_2^{-3})$. Clearly, $|\xi_1\overrightarrow{C}y|\le\overline{p}+1$. By Claim 2, $y\xi_2^+\not\in E$. Next, by Claim 4, $y\xi_2\not\in E$. Further, if $y\xi_2^-\in E$ then
$$
\xi_1x_1\overrightarrow{P}\xi_2\xi_2\xi_2^+\xi_2^-y\overrightarrow{C}\xi_2^{-2}\overrightarrow{Q}v_2\xi_1
$$
is longer than $C$, a contradiction. Finally, since $\mu(\Upsilon)=1$, we have $yv\not\in E$ for each $v\in V(\xi_2^{+2}\overrightarrow{C}\xi_1^{-})$. But then $G\backslash \{\xi_1,\xi_2^{-2}\}$ has at least three components, contradicting $\tau\ge1$.\\

 \textbf{Case a1.1.4}. $v_1=\xi_1$.

If $v_2v_3\in E$ for some $v_3\in V(\xi_2^{+2}\overrightarrow{C}\xi_1^{-})\cup V(\xi_1^+\overrightarrow{C}\xi_2^{-2})$ then we can argue as in Cases a1.1.1-a1.1.3. Otherwise $v_2v_3\in E$ for some $v_3\in\{\xi_2^-,\xi_2^+,\xi_2\}$. If $v_3\in\{\xi_2,\xi_2^+\}$ then we can show, as in Case a1.1.3, that $G\backslash \{\xi_1,v_3\}$ has at least three components, contradicting $\tau\ge1$.
     Now let $v_3=\xi_2^-$. Consider an arbitrary vertex $v\in V(Q)\backslash \{v_1\}$. Since $C$ is extreme, we have $N(v)\cap \{\xi_2,\xi_2^+\}=\emptyset$. Next, if $vy\in E$ for some $y\in V(C)\backslash \{\xi_1,\xi_2,\xi_2^-,\xi_2^+\}$ then we can argue as in Cases a1.1.1-a1.1.3. Thus, we can assume that $N(v)\subseteq V(Q)\cup \{\xi_2^-\}$, implying that $|Q|\ge\delta-1=\overline{p}+1$. Let $w\in V(\xi_1^+\overrightarrow{C}\xi_2^{-3})$. Since $|\xi_1\overrightarrow{C}w|\le \overline{p}+1$, we have $w\xi_2^+\not\in E$ (by Claim 2) and $w\xi_2\not\in E$ (by Claim 4). Recalling also that $\mu(\Upsilon)=1$, we conclude that $N(v)\subseteq V(\xi_1\overrightarrow{C}\xi_2^{-})$. If $\xi_2^{-2}\xi_2,\xi_2^{-2}\xi_2^+\not\in E$ then clearly $G\backslash \{\xi_1,\xi_2^-\}$ has at least three components, contradicting $\tau\ge1$. Hence, either $\xi_2^{-2}\xi_2\in E$ or $\xi_2^{-2}\xi_2^+\in E$. \\

\textbf{Case a1.1.4.1}. $\xi_2^{-2}\xi_2\in E$.

If $\xi_2^{-2}\xi_2^+\not\in E$ then $G\backslash \{\xi_1,\xi_2,\xi_2^-\}$ has at least four components, contradicting $\tau\ge1$. Hence, $\xi_2^{-2}\xi_2^+\in E$, that is $\langle\xi_2,\xi_2^-,\xi_2^{-2},\xi_2^+\rangle$ is a complete graph. If $V(G)=V(C\cup P\cup Q)$ then $n=4\delta+1\equiv 1 (mod\ 4)$ with $c=2\delta+3$, and we are done. Otherwise, as in previous cases, we can show that $\tau<1$, a contradiction.\\

\textbf{Case a1.1.4.2}. $\xi_2^{-2}\xi_2^+\in E$.

If $\xi_2^{-2}\xi_2\not\in E$ then $G\backslash \{\xi_1,\xi_2^-,\xi_2^+\}$ has at least four components, contradicting $\tau\ge1$. Otherwise $\langle\xi_2,\xi_2^-,\xi_2^{-2},\xi_2^+\rangle$ is a complete graph and we can argue as in Case a1.1.4.1. \\

\textbf{Case a1.1.5}. $v_1\in \{\xi_2,\xi_2^-,\xi_2^+\}$.

Since $C$ is extreme, we have $v_2\not\in \{\xi_2,\xi_2^-,\xi_2^+\}$ and therefore, we can argue as in Cases a1.1.1-1.1.4.  \\

\textbf{Case a1.2}.  $|I_1|=|I_2|=\overline{p}+3$.

We can show that $n=3\delta+1\equiv 1(mod\ 3)$ with $c=2\delta+2$, by arguing as in Case a1.1.\\

\textbf{Case a2}. $s=3$.

By Claim 6, $|I_1|=|I_2|=|I_3|=\overline{p}+3=\delta$ and $\xi_2^-\xi_2^+\in E$. If $\delta\ge4$ then $c=3\delta\ge2\delta+4$, contradicting (1). Hence $\delta=3$ and therefore, $\overline{p}=0$. Put
$$
C=\xi_1w_1w_2\xi_2w_3w_4\xi_3w_5w_6\xi_1,
$$
where $w_2w_3\in E$. Using Claims 2-5, we can show that 
It is not hard to see that
$$
N_C(w_1)=\{w_2,\xi_1,\xi_3\}, \ N_C(w_6)=\{w_5,\xi_1,\xi_3\}.
$$
Analogous relations hold for $w_4,w_5$. If $V(G\backslash C)=\{x_1\}$ then $n=10\equiv1(mod\ 3)$ with $c=9=2\delta+3>2\delta+2$, and we are done. Otherwise $N(y)=\{v_1,v_2,v_3\}$ for some $y\in V(G\backslash C)\backslash\{x_1\}$ with $N(y)\subseteq V(C)$. Since $C$ is extreme, it is not hard to see that either $N(y)=\{w_2,\xi_1,\xi_3\}$ or $N(y)=\{w_3,\xi_1,\xi_3\}$ or $N(y)=\{\xi_1,\xi_2,\xi_3\}$. But then $G\backslash N(y)$ has at least four components, contradicting $\tau\ge1$. Claim 8 is proved.\\

\textbf{Claim 9}. If $\mu=3$ then $G$ is the Petersen graph, that is $n=10\equiv1(mod\ 3)$ with $c\ge2\delta+2$. 

 \textbf{Proof}. By Claim 5, $\Upsilon(I_1,...,I_s)$ contains three pairwise crossing intermediate independent edges $e_1,e_2,e_3$. Let $e_1=w_1w_2$, $e_2=w_3w_4$ and $e_3=w_5w_6$. If $w_1,w_3,w_5\in V(I_f^\ast)$ for some $f\in \{1,...,s\}$, then we can argue as in proof of Claim 7. Otherwise we can assume w.l.o.g. that $w_1,w_3\in V(I_f^\ast)$, $w_2,w_5\in V(I_g^\ast)$ and $w_4,w_6\in V(I_h^\ast)$ for some distinct $f,g,h\in\{1,...,s\}$, where both $\xi_f,\xi_g,\xi_h$ and $w_1,w_3,w_5,w_2,w_4,w_6$ occur on $\overrightarrow{C}$ in a consecutive order. By Claim 1(a1 and a5), $|I_f|=|I_g|=|I_h|=\overline{p}+3$ and $|I_i|=\overline{p}+2$ for each $i\in\{1,...,s\}\backslash\{f,g,h\}$. Put
$$
|\xi_f\overrightarrow{C}w_1|=d_1, \  |w_1\overrightarrow{C}w_3|=d_2,  \  |w_3\overrightarrow{C}\xi_{f+1}|=d_3,
$$
$$
|\xi_g\overrightarrow{C}w_5|=d_4, \  |w_5\overrightarrow{C}w_2|=d_5,  \  |w_2\overrightarrow{C}\xi_{g+1}|=d_6,
$$
$$
|\xi_h\overrightarrow{C}w_4|=d_7, \  |w_4\overrightarrow{C}w_6|=d_8,  \  |w_6\overrightarrow{C}\xi_{h+1}|=d_9.
$$
If $d_3+d_7\ge\overline{p}+3$, $d_1+d_6\ge\overline{p}+3$ and $d_4+d_9\ge\overline{p}+3$ then clearly $|I_f|+|I_g|+|I_h|\ge3\overline{p}+12$, a contradiction. Otherwise we can assume w.l.o.g. that $d_3+d_7\le\overline{p}+2$. Further, if either $d_1\ge2$ or $d_9\ge2$ then we can argue as in the proof of Claim 7 (Case a1.1). Hence, we can assume that $d_1=d_9=1$. By Claim 2, $d_4=d_6=1$. By the same reason, using the fact that $d_1=d_6=1$, we get $d_3=d_7=1$. \\

\textbf{Case a1}. Either $\xi_{h+1}\not=\xi_f$ or $\xi_{f+1}\not=\xi_g$ or $\xi_{g+1}\not=\xi_h$.

Assume w.l.o.g. that $\xi_{h+1}\not=\xi_f$, implying that $|I_{f-1}|=\overline{p}+2$.  By Claim 5, $\xi_f^-y\not\in E$ for each $y\in V(I_i^\ast)$ and $i\in\{1,...,s\}\backslash\{f-1\}$. Moreover, by Claim 4, $\xi_f^-y\not\in E$ for each $y\in\{\xi_{f+1},\xi_h\}$. If $N(\xi_f^-)\subseteq V(C)$ then $d(\xi_f^-)\le\delta-1$, a contradiction. Otherwise we can argue as in the proof of Claim 6 (Case a1.2.1). \\

\textbf{Case a2}. $\xi_{h+1}=\xi_f$, $\xi_{f+1}=\xi_g$, $\xi_{g+1}=\xi_h$.

It follows that $s=3$. Assume w.l.o.g. that $f=1$, $g=2$ and $h=3$.\\

\textbf{Case a2.1.} Either $d_2\ge2$ or $d_5\ge2$ or $d_8\ge2$.

Assume w.l.o.g. that $d_2\ge2$, that is $w_1^+\not=w_3$. If  $\overline{p}=0$ then $|I_1|=3$, implying that $d_2=1$, a contradiction. Let $\overline{p}\ge1$. By Claim 4, $w_1^+\xi_2,w_1^+\xi_3\not\in E$. If $N(w_1^+)\subseteq V(C)$ then by Claim 4, $N(w_1^+)\subseteq V(w_1^{+2}\overrightarrow{C}w_3)\cup \{\xi_1\}$. Since $|I_1|=\overline{p}+3$, we have $|w_1^+\overrightarrow{C}w_3|\le \overline{p}$. But then $d(w_1^+)\le\overline{p}+1=\delta-2$, a contradiction. If $N(w_1^+)\not\subseteq V(C)$ then we can argue as in the proof of Claim 6 (Case a1.2.1). \\

\textbf{Case a2.2.} $d_2=d_5=d_8=1$.

It follows that $|I_i|=3$ $(i=1,2,3)$, that is $\overline{p}=0$, $\delta=3$ and $c=9$.   Clearly $\langle V(C)\cup \{x_1\}\rangle$ is the Petersen graph. If $V(G\backslash C)\not=\{x_1\}$ then it is not hard to see that $c\ge10$, a contradiction. Otherwise, $n=10\equiv 1(mod\ 3)$ with $c=9=2\delta+3>2\delta+2$. Claim 9 is proved.\\

Thus, the result holds from Claims 7,8,9.\\

\textbf{Case 2}. $\overline{p}=\delta-1$.

Clearly, $|N_C(x_i)|\ge1$ $(i=1,2)$.\\

\textbf{Case 2.1}. $x_1y_1, x_2y_2\in E$ for some distinct $y_1,y_2\in V(C)$.

We distinguish three main subcases.\\

\textbf{Case 2.1.1}. There exists a path $Q=z\overrightarrow{Q}y$ with $z\in V(P)$, $y\in V(C)\backslash \{y_1,y_2\}$ and $V(Q)\cap V(C\cup P)=\{z,y\}$.

Assume w.l.o.g. that $y\in V(y_1^+\overrightarrow{C}y_2^-)$. Since $C$ is extreme, we have 
$$
|y_1\overrightarrow{C}y|\ge|x_1\overrightarrow{P}z|+2, \  |y\overrightarrow{C}y_2|\ge|z\overrightarrow{P}x_2|+2,  \  |y_2\overrightarrow{C}y_1|\ge\delta+1.
$$
By summing, we get $|C|\ge2\delta+4$, contradicting (1).\\

\textbf{Case 2.1.2}. There exists a path $Q=z\overrightarrow{Q}y$ with $z\in V(y_1^+\overrightarrow{C}y_2^-)$, $y\in V(y_2^+\overrightarrow{C}y_1^-)$ and $V(Q)\cap V(C\cup P)=\{z,y\}$. 

By Claim 1(a1), $|C|\ge2\overline{p}+6=2\delta+4$, contradicting (1).\\

\textbf{Case 2.1.3}. $G\backslash \{y_1,y_2\}$ has at least three components.

It follows that $\tau<1$, contradicting the hypothesis.\\

\textbf{Case 2.2}. $N_C(x_1)=N_C(x_2)=\{y\}$ for some $y\in V(C)$.

It follows that
$$
N(x_1)=(V(P)\cup\{y\})\backslash \{x_1\},  \  N(x_2)=(V(P)\cup\{y\})\backslash \{x_2\}.
$$
Moreover, $x_1\overrightarrow{P}v^-x_2\overleftarrow{P}v$ is a longest path in $G\backslash C$ for each $v\in V(x_1^+\overrightarrow{P}x_2)$. Since $G$ is 2-connected, we have $wz\in E$ for some $w\in V(P)$ and $z\in V(C)\backslash\{y\}$. If $w=x_1$ then using the path $zx_1\overrightarrow{P}x_2y$, we can argue as in Case 2.1. Otherwise we can use the path $yx_1\overrightarrow{P}w^-x_2\overleftarrow{P}wz$.\\

\textbf{Case 3}. $\overline{p}\ge\delta$.

\textbf{Case 3.1}. $x_1y_1,x_2y_2\in E$ for some distinct $y_1,y_2\in V(C)$.

Clearly, $|y_1\overrightarrow{C}y_2|\ge\delta+2$ and $|y_2\overrightarrow{C}y_1|\ge\delta+2$, which yields $|C|\ge2\delta+4$, contradicting (1).\\

\textbf{Case 3.2}. $N_C(x_1)=N_C(x_2)=\{y\}$ for some $y\in V(C)$.

Let $y_1,y_2,...,y_t$ be the elements of $N_P^+(x_2)$ 
occuring on $\overrightarrow{P}$ in a consecutive order.
 Put $H=\langle V(y_1^-\overrightarrow{P}x_2)\rangle$ and
$$
P_i=x_1\overrightarrow{P}y_i^-x_2\overleftarrow{P}y_i \ (i=1,...,t).
$$
Since $P_i$ is a longest path in $G\backslash C$ for each $i\in\{1,...,t\}$, we can assume w.l.o.g. that $P$ is chosen so as to maximize $|V(H)|$. If $y_iz\in E$ for some $i\in\{1,...,t\}$ and $z\in V(C)\backslash\{y\}$, then we can argue as in Case 3.1. Otherwise $N(y_i)\subseteq V(H)\cup \{y\}$ $(i=1,...,t)$, that is $|N_H(y_i)|\ge\delta-1$ $(i=1,...,t)$. By Lemma 3, for each distinct $u,v\in V(H)$, there is a path in $H$ of length at least $\delta-1$, connecting $u$ and $v$. Since $G$ is 2-connected, $H$ and $C$ are connected by two vertex disjoint paths. This means that there is a path $Q=y_1\overrightarrow{Q}y_2$ of length at least $\delta+1$ with $V(Q)\cap V(C)=\{y_1,y_2\}$. Further, we can argue as in Case 2.1. \\

\textbf{Case 3.3}. Either $N_C(x_1)=\emptyset$ or $N_C(x_2)=\emptyset$.

Assume w.l.o.g. that $N_C(x_1)=\emptyset$. By arguing as in Case 3.2, we can find a path $Q=y_1\overrightarrow{Q}y_2$ of length at least $\delta+2$ with $V(Q)\cap V(C)=\{y_1,y_2\}$, and the result follows immediately. Theorem 1 is proved.        \quad   \quad         \rule{7pt}{6pt}

\noindent Institute for Informatics and Automation Problems\\ National Academy of Sciences\\
P. Sevak 1, Yerevan 0014, Armenia\\ E-mail: zhora@ipia.sci.am


\begin{thebibliography}{10}

\bibitem{[1]}   D. Bauer and E. Schmeichel, Long cycles in tough graphs, Technical Report 8612, Stevens Institute of Technology, 1986. 

\bibitem{[2]}   J.A. Bondy and U.S.R. Murty, Graph Theory with Applications. Macmillan, London and Elsevier, New York (1976).

\bibitem{[3]}  V. Chv\'{a}tal, Tough graphs and hamiltonian circuits, Discrete Math. 5(1973) 215-228.

\bibitem{[4]}   G. A. Dirac, Some theorems on abstract graphs, Proc. London, Math. Soc., 2 (1952) 69-81.

\bibitem{[5]}   H.A. Jung, On maximal circuits in finite graphs, Ann. Discrete Math. 3 (1987) 129-144.

\bibitem{[6]}   H.-J. Voss, Bridges of longest circuits and of longest paths in graphs, Beitrage zur Graphentheorie und deren Anwendungen, Vorgetr. auf dem. int. Kolloq., Oberhof (DDR) (1977) 275-286. 

\end{thebibliography}
\end{document}